\input amstex
\documentstyle{amsppt}
\document
\topmatter
\title
K\"ahler surfaces with quasi-constant  holomorphic curvature.
\endtitle
\author
W\l odzimierz Jelonek
\endauthor

\abstract{The aim of this paper is to describe   K\"ahler surfaces
with quasi-constant holomorphic sectional curvature. }
\thanks{MS Classification: 53C55,53C25,53B35. Key words and phrases:
K\"ahler surface, holomorphic sectional curvature, quasi constant
holomorphic sectional curvature,QCH manifold, ambik\"ahler
mani\-fold }\endthanks
 \endabstract
\endtopmatter
\define\G{\Gamma}
\define\bJ{\overline{J}}
\define\DE{\Cal D^{\perp}}

\define\n{\nabla}
\define\om{\omega}

\define\r{\rightarrow}
\define\w{\wedge}
\define\k{\diamondsuit}
\define\th{\theta}
\define\p{\partial}
\define\a{\alpha}

\define\lb{\lambda}

\define\1{D_{\lb}}
\define\2{D_{\mu}}
\define\0{\Omega}

\define\De{\Cal D}

\define\m{(M,g,J)}
\define \E{\Cal E}
\bigskip
{\bf 0. Introduction. } The aim of the present paper is to
describe connected K\"ahler surfaces $(M,g,J)$ admitting a global,
$2$-dimensional, $J$-invariant  distribution $\De$ having the
following property:  The holomorphic curvature
$K(\pi)=R(X,JX,JX,X)$ of any $J$-invariant $2$-plane $\pi\subset
T_xM$, where $X\in \pi$ and $g(X,X)=1$, depends only on the point
$x$ and the number $|X_{\De}|=\sqrt{g(X_{\De},X_{\De})}$, where
$X_{\De}$ is an orthogonal projection of $X$ on $\De$. In this
case  we have
$$R(X,JX,JX,X)=\phi(x,|X_{\De}|)$$ where $\phi(x,t)=a(x)+b(x)t^2+c(x)t^4$ and
 $a,b,c$ are smooth functions on $M$. Also $R=a\Pi+b\Phi+c\Psi$
 for certain curvature tensors $\Pi,\Phi,\Psi\in \bigotimes^4\frak X^*(M)$
  of K\"ahler type. The investigation of such manifolds, called QCH K\"ahler manifolds,  was
started by G. Ganchev and V. Mihova in [G-M-1],[G-M-2]. In our
paper [J-2] we  used their local results to obtain a global
classification of such manifolds under the assumption that $\dim
M=2n\ge 6$.
 By $\E$ we shall denote the $2$-dimensional distribution
 which is the orthogonal complement of $\De$ in $TM$.    In
 the present paper we show that a K\"ahler surface $\m$ is a QCH
 manifold with respect to a distribution $\De$ if and only if is a QCH manifold with respect to the distribution $\E$.
 We also prove that $\m$ is a QCH K\"ahler surface if and only if the antiselfdual Weyl tensor $W^-$ is degenerate and
 there exist a negative almost complex
 structure $\overline{J}$ which preserves the Ricci tensor $Ric$
 of $\m$ i.e. $Ric(\bJ.,\bJ.)=Ric(.,.)$ and such that
 $\overline{\om}=g(\overline{J}.,.)$ is an eigenvector of $W^-$
 corresponding to simple eigenvalue of $W^-$.  Equivalently
 $(M,g,J)$ is a QCH K\"ahler surface iff it admits a negative
 almost complex structure $ \overline{J}$ satisfying the Gray
 second condition
 $R(X,Y,Z,W)-R(\overline{J}X,\overline{J}Y,Z,W)=R(\overline{J}X,Y,\overline{J}Z,W)+R(\overline{J}X,Y,Z,\overline{J}W)$.
  In [A-C-G-1]  Apostolov, Calderbank and Gauduchon have classified weakly selfdual K\"ahler surfaces, extending
  the result of Bryant who classified self-dual K\"ahler surfaces [B].
  Weakly self-dual K\"ahler surfaces turned out to be of Calabi
  type and of orthotoric type or surfaces with parallel Ricci tensor.
   We show that any Calabi type K\"ahler surface and every
 orthotoric K\"ahler surface is a QCH manifold. In both cases the
 opposite complex strucure $\bJ$ is conformally K\"ahler. We also
 classify locally homogeneous QCH K\"ahler surfaces.
\bigskip
{\bf 1.   Almost complex structure $ \overline{J}$.} Let $\m$ be a
$4$-dimensional K\"ahler manifold with a $2$-dimensional
$J$-invariant distribution $\De$. Let $\frak X(M)$ denote the
algebra of all differentiable vector fields on $M$ and $\G(\De)$
denote the set of local sections of the distribution $\De$. If
$X\in\frak X(M)$ then by $X^{\flat}$ we shall denote the 1-form
$\phi\in \frak X^*(M)$ dual to $X$ with respect to $g$, i.e.
$\phi(Y)=X^{\flat}(Y)=g(X,Y)$. By $\om$ we shall denote the
K\"ahler form of $\m$ i.e. $\om(X,Y)=g(JX,Y)$. Let $\m$ be a QCH
K\"ahler surface with respect to $J-invariant$ 2-dimensional
distribution $\De$. Let us denote by $\E$ the distribution $\DE$,
which is a $2$-dimensional, $J$-invariant distribution. By $h,m$
respectively we shall denote the tensors $h=g\circ (p_{\De}\times
p_{\De}),m=g\circ (p_{\E}\times p_{\E})$, where $p_{\De},p_{\E}$
are the orthogonal projections on $\De,\E$ respectively. It
follows that $g=h+m$. Let us define almost complex structure $
\overline{J}$ by $ \overline{J}_{|\E}=-J_{|\E}$ and
$\overline{J}_{|\De}=J_{|\De}$. Let $\th(X)=g(\xi,X)$ and
$J\th=-\th\circ J$ which means that $J\th(X)=g(J\xi,X)$. For every
almost Hermitian manifold $\m$ the self-dual Weyl tensor $W^+$
decomposes under the action of the unitary group $U(2)$. We have
$\bigwedge^*M=\Bbb R\oplus LM$ where $LM=[[\bigwedge^{(0,2)}M]]$
and we can write $W^+$ as a matrix with respect to this block
decomposition

$$W^+=\pmatrix \frac{\kappa}6& W_2^+\cr (W_2^+)^*&W_3^+-\frac{\kappa}{12}Id_{|LM}\endpmatrix$$
where $\kappa$ is the conformal scalar curvature of $(M,g,J)$ (see
[A-A-D]). The selfdual Weyl tensor $W^+$ of $\m$ is called
degenerate if $W_2=0,W_3=0$. In general the self-dual Weyl tensor
of 4-manifold $(M,g)$ is called degenerate if it has at most two
eigenvalues as an endomorphism $W^+:\bigwedge^+M\r\bigwedge^+M$.
We say that an almost Hermitian structure $J$ satisfies the second
Gray curvature condition if
$$R(X,Y,Z,W)-R(JX,JY,Z,W)=R(JX,Y,JZ,W)+R(JX,Y,Z,JW),\tag {G2}$$
which is equivalent to $Ric(J,J)=Ric$ and $W_2^+=W_3^+=0$.  Hence
$\m$ satisfies the second Gray condition if $J$ preserves the
Ricci tensor and $W^+$ is degenerate.  We shall denote by $Ric_0$
and $\rho_0$ the trace free part of the Ricci tensor $Ric$ and the
Ricci form $\rho$ respectively.   An ambik\"ahler structure on a
real 4-manifold consists of a pair of K\"ahler metrics
$(g_+,J_+,\om_+)$ and $(g_-,J_-,\om_-)$ such that $g_+$ and $g_-$
are conformal metrics and $J_+$ gives an opposite orientation to
that given by $J_-$ (i.e the volume elements $\frac12\om_+\w\om_+$
and $\frac12\om_-\w\om_-$ have opposite signs).

\medskip
{\bf 2. Curvature tensor of a QCH K\"ahler surface.} We shall
recall some results from [G-M-1]. Let
$$R(X,Y)Z=([\n_X,\n_Y]-\n_{[X,Y]})Z\tag 2.1$$ and let us write $$R(X
,Y,Z,W)=g(R(X,Y)Z,W).$$ If $R$ is the curvature tensor of a QCH
K\"ahler manifold $\m$, then there exist functions $a,b,c\in
C^{\infty}(M)$ such that
$$R=a\Pi+b\Phi+c\Psi,\tag 2.2$$
where $\Pi$ is the standard K\"ahler tensor of constant
holomorphic curvature i.e.
$$\gather \Pi(X,Y,Z,U)=\frac14(g(Y,Z)g(X,U)-g(X,Z)g(Y,U)\tag 2.3\\
+g(JY,Z)g(JX,U)-g(JX,Z)g(JY,U)-2g(JX,Y)g(JZ,U)),\endgather $$
the tensor $\Phi$ is defined by the following relation
$$\gather \Phi(X,Y,Z,U)=\frac18(g(Y,Z)h(X,U)-g(X,Z)h(Y,U)\tag 2.4\\+g(X,U)h(Y,Z)-g(Y,U)h(X,Z)
+g(JY,Z)h(JX,U)\\-g(JX,Z)h(JY,U)+g(JX,U)h(JY,Z)-g(JY,U)h(JX,Z)\\
-2g(JX,Y)h(JZ,U)-2g(JZ,U)h(JX,Y)),\endgather$$ and finally
$$\Psi(X,Y,Z,U)=-h(JX,Y)h(JZ,U)=-(h_J\otimes h_J)(X,Y,Z,U).\tag 2.5$$
where $h_J(X,Y)=h(JX,Y)$.
 Let $V=
 (V,g,J)$ be a real $2n$ dimensional vector space with
complex structure $J$ which is skew-symmetric with respect to the
scalar product $g$ on $V$.  Let assume further that  $V= D\oplus
E$ where $D$ is a 2-dimensional, $J$-invariant subspace of $V$,
$E$ denotes its orthogonal complement in $V$. Note that  the
tensors $\Pi,\Phi,\Psi$ given above are of K\"ahler type.  It is
easy to check that for a unit vector $X\in V$
$\Pi(X,JX,JX,X)=1,\Phi(X,JX,JX,X)=|X_{D}|^2,\Psi(X,JX,JX,X)=|X_{D}|^4$,
where $X_D$ means an orthogonal projection of a vector $X$ on the
subspace $D$ and $|X|=\sqrt{g(X,X)}$. It follows that for a tensor
$(2.2)$ defined on $V$ we have
$$R(X,JX,JX,X)=\phi(|X_D|)$$ where $\phi(t)=a+bt^2+ct^4$.

Let $J, \overline{J}$ be hermitian, opposite orthogonal structures
on a Riemannian 4-mani\-fold $(M,g)$ such that $J$ is a positive
almost complex structure. Let $\Cal
E=ker(J\overline{J}-Id),\De=ker(J\overline{J}+Id)$ and let the
tensors $\Pi,\Phi,\Psi$ be defined as above where
$h=g(p_{\De},p_{\De})$. Let us define a tensor
$K=\frac16\Pi-\Phi+\Psi$.  Then $K$ is a curvature tensor,
$b(K)=0,c(K)=0$ where $b$ is Bianchi operator and $c$ is the Ricci
contraction. Define the endomorphism
$K:\bigwedge^2M\r\bigwedge^2M$ by the formula
$g(K\phi,\psi)=-K(\phi,\psi)$ (see (2.1)). Then we have

\medskip
{\bf Lemma 1.} {\it The tensor $K$ satisfies $K(\bigwedge^+M)=0$.
Let $\phi,\psi\in\bigwedge^-M$ be the local forms  orthogonal to $
\overline{\om}$ such that $g(\phi,\psi)=g(\psi,\psi)=2$ and
$g(\phi,\psi)=0$.  Then $K( \overline{\om})=\frac13
\overline{\om}, K(\phi)=-\frac16\phi,K(\psi)=-\frac16\psi$.}
\medskip
{\it Proof.} A straightforward computation.$\k$
\medskip
In the special case of a K\"ahler surface $\m$ we get for a QCH
manifold $\m$
\medskip
{\bf Proposition  1.} {\it  Let $\m$ be a K\"ahler surface which
is a QCH manifold with respect to the distribution $\De$.  Then
$\m$ is also QCH manifold with respect to the distribution $\Cal
E=\De^{\perp}$ and if $\Phi',\Psi'$ are the above tensors with
respect to $\Cal E$ then}
$$R=(a+b+c)\Pi-(b+2c)\Phi'+c\Psi'.\tag 2.6$$

\medskip
{\it Proof.}  Let us assume that $X\in TM,||X||=1$.  Then if
$\a=||X_{\De}||,\beta=||X_{\E}||$ then $1=\a^2+\beta^2$. Hence
$R(X,JX,JX,X)=a+b\a^2+c\a^4=a+b(1-\beta^2)+c(1-\beta^2)^2=a+b+c-(b+2c)\beta^2+c\beta^4$.$\k$

\medskip
 If $\m$ is a QCH K\"ahler surface then one can show that
the Ricci tensor $\rho$ of $\m$ satisfies the equation
$$\rho(X,Y)=\lb m(X,Y)+\mu h(X,Y)\tag 2.7$$
where $\lb=\frac{3}2a+\frac  b4,\mu=\frac{3}2a+\frac{5}4b+c$ are
eigenvalues of $\rho$ (see [G-M-1], Corollary 2.1 and Remark 2.1.)
In particular the distributions $\E,\De$ are eigendistributions of
the tensor $\rho$ corresponding to the eigenvalues $\lb,\mu$ of
$\rho$. The Kulkarni-Nomizu product of two symmetric
$(2,0)$-tensors $h,k\in \bigotimes^2TM^*$ we call a tensor
$h\oslash k$ defined as follows:

$$\gather h\oslash
k(X,Y,Z,T)=h(X,Z)k(Y,T)+h(Y,T)k(X,Z)\\-h(X,T)k(Y,Z)-h(Y,Z)k(X,T).\endgather$$
Similarly we define the Kulkarni-Nomizu product of two 2-forms
$\om,\eta$

$$\gather \om\oslash
\eta(X,Y,Z,T)=\om(X,Z)\eta(Y,T)+\om(Y,T)\eta(X,Z)\\-\om(X,T)\eta(Y,Z)-\om(Y,Z)\eta(X,T).\endgather$$
Then $b(\om\oslash\eta)=-\frac23\om\w\eta$ where $b$ is the
Bianchi operator. In fact

$$\gather 3b(\om\oslash
\eta)(X,Y,Z,T)=\om(X,Z)\eta(Y,T)+\om(Y,T)\eta(X,Z)-\om(X,T)\eta(Y,Z)\\-\om(Y,Z)\eta(X,T)
+\om(Y,X)\eta(Z,T)+\om(Z,T)\eta(Y,X)\\-\om(Y,T)\eta(Z,X)-\om(Z,X)\eta(Y,T)
+\om(Z,Y)\eta(X,T)\\+\om(X,T)\eta(Z,Y)-\om(Z,T)\eta(X,Y)-\om(X,Y)\eta(Z,T)\\=-2\om\w\eta(X,Y,Z,T)
.\endgather$$

 Note that
$$\gather \Pi=-\frac14(\frac12(g\oslash
g+\om\oslash\om)+2\om\otimes\om)),\tag 2.8\\
\Phi=-\frac18(h\oslash g + h_J\oslash \om +2\om\otimes
h_J+2h_J\otimes
\om),\tag 2.9\\
\Psi=-h_J\otimes h_J,\tag 2.10\endgather$$ where $\om=g(J.,.)$ is
the K\"ahler form. Note that $b(\Psi)=\frac13h_J\w h_J=0$ since
$h_J=e_1\w e_2$ is primitive, where $e_1,e_2$ is an orthonormal
basis in $\De$.

\medskip
{\bf Theorem  1. } {\it Let $\m$ be a K\"ahler surface.  If $\m$
is a QCH manifold then $W^-=c(\frac16\Pi-\Phi+\Psi)$ and $W^-$ is
degenerate.  The 2-form $\overline{\om}$ is an eigenvector of
$W^-$ corresponding to a simple eigenvalue of $W^-$ and $
\overline{J}$ preserves the Ricci tensor. On the other hand let us
assume that $\m$ admits a negative almost complex structure $
\overline{J}$ such that $Ric(\overline{J},\overline{J})=Ric$. Let
$\Cal E=ker(J\overline{J}-Id),\De=ker(J\overline{J}+Id)$. If
$W^-=\frac{\kappa}2(\frac16\Pi-\Phi+\Psi)$ or equivalently if the
half-Weyl tensor $W^-$ is degenerate and $ \overline{\om}$ is an
eigenvector of $W^-$ corresponding to a simple eigenvalue of $W^-$
then   $\m$ is a QCH manifold.}

{\it Proof.} Note that for a K\"ahler surface $\m$ the Bochner
tensor coincides with $W^-$ and we have

$$\gather  R=-\frac{\tau}{12}(\frac14(g\oslash
g+\om\oslash\om)+\om\otimes\om)\\-\frac1{4}(\frac12(Ric_0\oslash g
+\rho_0\oslash\om)+\rho_0\otimes\om+\om\otimes\rho_0)+W^-.\endgather
$$

If $\m$ is a QCH K\"ahler surface then $Ric=\lb m+ \mu h$ where
$\lb=\frac32a+\frac b4,\mu= \frac32a+\frac 54b+c$.  Consequently
$Ric_0=-\frac{b+c}2m+\frac{b+c}2h=\delta h-\delta m$ where
$\delta=\frac{b+c}2$. Hence $Ric_0=2\delta h-\delta g$. Hence we
have

$$\gather  R=-\frac{\tau}{12}(\frac14(g\oslash
g+\om\oslash\om)+\om\otimes\om)\\-\frac1{4}(\frac12((2\delta
h-\delta g)\oslash g +(2\delta h_J-\delta \om)\oslash\om)+(2\delta
h_J-\delta \om)\otimes\om+\\\om\otimes(2\delta h_J-\delta
\om))+W^-.\endgather
$$
Consequently
$$R=\frac{\tau}6\Pi+2\delta\Phi-\delta\Pi+W^-=(a-\frac c6)\Pi+(b+c)\Phi+W^-$$
and $a\Pi+b\Phi+c\Psi=(a-\frac c6)\Pi+(b+c)\Phi+W^-$ hence
$W^-=c(\frac16\Pi-\Phi+\Psi)$.   It follows that $W^-$ is
degenerate and $ \overline{\om}$ is an eigenvalue of $W^-$
corresponding to the simple eigenvalue of $W^-$.  It is also clear
that $Ric( \overline{J},\overline{J})=Ric$.

On the other hand let us assume that a K\"ahler surface $\m$
admits a negative almost complex structure $ \overline{J}$
preserving the Ricci tenor $Ric$ and such that $W^-$ is degenerate
with eigenvector $ \overline{\om}$ corresponding to the simple
eigenvalue of $W^-$.   Equivalently it means that $\overline{J}$
satisfies the second Gray condition of the curvature i.e.
$R(X,Y,Z,W)-R(\overline{J}X,\overline{J}Y,Z,W)=R(\overline{J}X,Y,\overline{J}Z,W)+R(\overline{J}X,Y,Z,\overline{J}W).$
Then  $W^-=\frac{\kappa}2((\frac16\Pi-\Phi+\Psi)$. If
$Ric_0=\delta(h-m)$ then as above
$R=\frac{\tau}6\Pi+2\delta\Phi-\delta\Pi+W^-$. Consequently
$R=(\frac{\tau}6-\delta)\Pi+2\delta\Phi+\frac{\kappa}2(\frac16\Pi-\Phi+\Psi)$
and consequently
$$R=(\frac{\tau}6-\delta+\frac\kappa{12})\Pi+(2\delta-\frac{\kappa}2)\Phi+\frac{\kappa}2\Psi.\tag 2.11$$
$\k$
\medskip
{\it Remark.}  Note that $\kappa$ is the conformal scalar
curvature of $(M,g, \overline{J})$.  The Bochner tensor of QCH
manifold was first identified in [G-M-2].

\medskip
{\it Corollary.}   A K\"ahler surface $\m$ is a QCH manifold iff
it admits a negative almost complex structure $\overline{J}$
satisfying the second Gray condition of the curvature i.e.
$$R(X,Y,Z,W)-R(\overline{J}X,\overline{J}Y,Z,W)=R(\overline{J}X,Y,\overline{J}Z,W)+R(\overline{J}X,Y,Z,\overline{J}W)$$
The $J$-invariant distribution $\De$ with respect to which $\m$ is
a QCH manifold is given by $\De=ker(J\overline{J}-Id)$ or by
$\De=ker(J\overline{J}+Id)$.

\medskip
{\bf Theorem 2. }{\it Let us assume that $\m$ is a K\"ahler
surface admitting a negative Hermitian structure $\overline{J}$
such that $Ric(\overline{J},\overline{J})=Ric$.  Then $\m$ is a
QCH manifold.}
\medskip
{\it Proof.}  If a Hermitian manifold $\m$ has a J-invariant Ricci
tensor $Ric$ then the tensor $W^+$ is degenerate (see [A-G]). $\k$

\medskip
{\it Remark.} If a K\"ahler surface $\m$ is compact and admits a
negative Hermitian structure $\overline{J}$ as above then
$(M,g,\overline{J})$ is locally conformally K\"ahler and hence
globally conformally K\"ahler if $b_1(M)$ is even. Thus  $\m$ is
ambiK\"ahler since $b_1(M)$ is even.
\medskip

Now we give examples of QCH K\"ahler surfaces. First we give (see
[A-C-G-1])
\medskip
{\it Definition.}  A K\"ahler surface $\m$ is said to be of Calabi
type if it admits a non-vanishing Hamiltonian  Killing vector
field $\xi$ such that the almost Hermitian pair $(g,I)$ -with $I$
equal to $J$ on the distribution spanned by $\xi$ and $J\xi$ and
$-J$ on the orthogonal distribution - is conformally K\"ahler.
\medskip

Every K\"ahler surface of Calabi type is given locally by

$$\gather
g=(az-b)g_{\Sigma}+w(z)dz^2+w(z)^{-1}(dt+\a)^2,\tag 2.12
\\\om=(az-b)\om_{\Sigma}+dz\w (dt+\a),
d\a=a\om_{\Sigma}\endgather$$ where $\xi=\frac{\p}{\p t}$.

The K\"ahler form of Hermitian structure $I$ is given by
$\om_I=(az-b)\om_{\Sigma}-dz\w (dt+\a)$ and the K\"ahler metric
corresponding to $I$ is $g_{-}=(az-b)^2g$.

If $a\ne0$ then the metric (*) is a product metric. If $a\ne0$
then we set $a=1,b=0$ and write $w(z)=\frac z{V(z)}$ hence
$$\gather
g=zg_{\Sigma}+\frac
z{V(z)}zdz^2+\frac{V(z)}z(dt+\a)^2,\tag2.13\\\om=z\om_{\Sigma}+dz\w
(dt+\a), d\a=\om_{\Sigma}\endgather$$

It is known that for a K\"ahler surface of Calabi type of
non-product type we have $\rho_0=\delta\om_I$  where
$\delta=-\frac1{4z}(\tau_{\Sigma}+(\frac{V_z}{z^2})_zz^2)$ (see
[A-C-G-1]) and consequently $Ric(I,I)=Ric$. This last relation
remains true in the product case metric.  Hence we have
\medskip

{\bf Theorem  3.  }{\it   Every K\"ahler surface of Calabi type is
a QCH K\"ahler surface.}

  \medskip
{\it  Definition.} A K\"ahler surface $\m$ is ortho-toric if it
admits two independent Hamiltonian Killing vector fields with
Poisson commuting momentum maps $\xi\eta$ and $\xi+\eta$ such that
$d\xi$ and $d\eta$ are orthogonal.

An explicit classification of ortho-toric K\"ahler metrics is
given in [A-C-G-1].  We have (this Proposition is proved in
[A-C-G-1], Prop.8 )
\medskip
{\bf Proposition.} {\it The almost Hermitian structure $(g,J,\om)$
defined by
$$\gather g=(\xi-\eta)(\frac{d\xi^2}{F(\xi)}-\frac{d\eta^2}{G(\eta)})+\frac1{\xi-\eta}(F(\xi)(dt+\eta dz)^2-G(\eta)(dt+\xi dz)^2\tag 2.14\\
Jd\xi=\frac{F(\xi)}{\xi-\eta}(dt+\eta dz),  Jdt=-\frac{\xi d\xi}{F(\xi)}-\frac{\eta d\eta}{G(\eta)}\tag 2.15\\
Jd\eta=-\frac{G(\eta)}{\xi-\eta}(dt+\xi dz),  Jdz=\frac{d\xi}{F(\xi)}+\frac{ d\eta}{G(\eta)},\\
\om=d\xi\w (dt+\eta dz)+d\eta\w(dt+\xi dz)\tag 2.16\endgather$$ is
orthotoric where $F,G$ are any functions of one variable. Every
orthotoric K\"ahler surface $\m$ is of this form.}

Any orthotoric surface has a negative Hermitian structure $
\overline{J}$, whose K\"ahler form $ \overline{\om}$ is given by
$$  \overline{\om}=d\xi\w (dt+\eta dz)-d\eta\w(dt+\xi dz)$$
and
$$\gather \overline{J}d\xi =Jd\xi=\frac{F(\xi)}{\xi-\eta}(dt+\eta dz),  \overline{J}dt=-\frac{\xi d\xi}{F(\xi)}+\frac{\eta d\eta}{G(\eta)}\tag 2.17\\
\overline{J}d\eta=Jd\eta=-\frac{G(\eta)}{\xi-\eta}(dt+\xi dz),
\overline{J}dz=\frac{d\xi}{F(\xi)}-\frac{
d\eta}{G(\eta)},\endgather $$

The structure $(g_-=(\xi-\eta)^2g,\overline{J})$ is K\"ahler.  We
also have $\rho_0=\delta \overline{\om}$ where
$\delta=\frac{F'(\xi)-G'(\eta)}{(2(\xi-\eta)^2}-\frac{F''(\xi)+G''(\eta)}{(4(\xi-\eta)}$.

In particular the Hermitian structure $ \overline{J}$ preserves
Ricci tensor $Ric$.  Hence we get

\medskip

{\bf Theorem  4.  }{\it Every orthotoric  K\"ahler surface is a
QCH K\"ahler surface.}
\medskip
Note that both Calabi type and orthotoric K\"ahler surfaces are
ambik\"ahler.  On the other hand we have
\bigskip
{\bf Theorem  5.  }  {\it Let $\m$ be ambi-K\"ahler surface which
is a QCH manifold. Then locally $\m$ is orthotoric or of Calabi
type or a product of two Riemannian surfaces or is an
anti-selfdual Einstein-K\"ahler surface.}
\medskip
{\it Proof.} (We follow [A-C-G-2]). Let us denote by $g_-$ the
second K\"ahler metric. Let us assume that $g_-\ne g$. Then
$g=\phi^{-2}g_-$ and the field $X=grad_{\om_-}\phi$ is a
 Killing vector field $L_Xg=L_Xg_-=0$ and is holomorphic  with respect to $
 \overline{J})$.     We shall show that $X$ is also holomorphic
 with respect to $J$.  In fact $Ric_0=\delta g(J \overline{J},.)$
 and $L_XRic=0, L_X\delta=0$.  Hence $0=\delta g((L_XJ)
 \overline{J},.)$ and consequently $L_XJ=0$ in
 $U=\{x:Ric_0(x)\ne0\}$.   If $(M,g)$  is Einstein then $W^+\ne0$
 everywhere or $(M,g,J)$ is anti-selfdual.  In the first case
 $X$ preserves the simple eigenspace  of $W^+$ and hence $\om$,
 cosequently
 $L_XJ=0$.

  Note that
 $X=\overline{J}grad_g\psi$ where $\psi=-\frac1{\phi}$. Since
 $L_X\om=0$ we have $d X\lrcorner\om=0$ and consequently the
 1-form $J\overline{J}d\psi$ is closed and locally equals
 $\frac12d\sigma$.
  Thus the two form
 $\0=\frac32\sigma\om+\psi^3\om_-$, where $\om_-$ is the K\"ahler
 form of $(M,g_-,\overline{J})$, is a Hamiltonian form in the sense of [A-C-G-1] and the
 result follows from the classification in [A-C-G-1]. This form is defined globally if $H^1(M)=0$.$\k$
\medskip
 {\it Remark.}  Note that in the compact case every Killing vector
 field on a K\"ahler surface is holomorphic.  If $\m$ is an
 Einstein
 K\"ahler anti-selfdual then in the case where it is not
 conformally flat the manifold $(M,g,\bJ)$ is a self-dual
 Einstein Hermitian conformal to self-dual K\"ahler metric. Such a
 metric must be either orthotoric or of Calabi type. Thus
  $\m$ is of Calabi type if $(M,g,\bJ)$ is of Calabi type, however
  $\m$ can not be orthotoric if $(M,g,\bJ)$ is orthotoric.

\medskip
 Now we shall investigate Einstein QCH K\"ahler surfaces.
\medskip
{\bf Theorem  6.}  {\it Let $\m$ be a K\"ahler-Einstein surface.
Then $\m$ is a QCH K\"ahler surface if and only if it admits a
negative Hermitian structure $ \overline{J}$ or it has constant
holomorphic curvature and admits any negative almost complex
structure. If $\m$ is QCH and the second case does not hold then
$\overline{J}$ is conformally K\"ahler hence $\m$ is ambiK\"ahler.
}
\medskip
{\it Proof.}  If an Einstein 4-manifold $(M,g)$ admits a
degenerate tensor $W^-$ then $W^-=0$ or $W^-\ne0$ on the whole of
$M$. In the second case by the result of Derdzinski it admits a
Hermitian structure $ \overline{J}$ which is conformally K\"ahler
and the metric $ (g(W^-,W^-))^{\frac13}g$ is a K\"ahler metric
with respect to $\overline{J}$.
\medskip
{\it Remark.} (Compare [A-C-G-1]). If $\m$ is a QCH K\"ahler
Einstein surface which is not anti-self-dual then in the case
$H^1(M)=0$ on $(M,g,J)$ there is defined global Hamiltonian two
form  surface and on the open and dense subset $U$ of $M$ the
metric $g$ is:

(a) a K\"ahler product metric of two Riemannian surfaces of the
same Gauss curvature

(b) K\"ahler Einstein metric of Calabi type over a Riemannian
surface $(\Sigma,g_{\Sigma})$ of constant scalar curvature $k$ of
the form (2.13) where $V(z)=a_1z^3+kz^2+a_2$

(c)  K\"ahler-Einstein ambitoric metric of parabolic type (see
[A-C-G-2])
\medskip
{\bf Theorem  7.}  {\it Let $\m$ be a self-dual K\"ahler surface
with $Ric_0\ne0$ everywhere on $M$. Then $\m$ is a QCH K\"ahler
surface with Hermitian complex structure $ \overline{J}$. }

{\it Proof.}  We show as in Th.1 that
$R=\frac{\tau}6\Pi+2\delta\Phi-\delta\Pi$ where $\rho_0=\delta
\overline{\om}$.  Note that in $U=\{x:Ric_0\ne0\}$ the negative
structure $ \overline{J}$ is uniquely determined and is Hermitian
in $U$ (see Prop.4 in [A-G]). $\k$
\medskip
{\it Remark.}  Note  that a selfdual K\"ahler surface $\m$ is QCH
if admits any negative almost complex structure $ \overline{J}$
preserving the Ricci tensor $Ric$. For example $\Bbb{CP}^2$ with
standard Fubini-Studi metric is selfdual however is not QCH since
it does not admit any negative almost complex structure.  However
the manifold $M=\Bbb{CP}^2-\{p_0\}$ for any point
$p_0\in\Bbb{CP}^2$ is QCH and admits a negative Hermitian complex
structure (see [J-3]). In [D-2] there are constructed many
examples of self-dual K\"ahler surfaces with $Ric_0\ne0$ hence QCH
K\"ahler self-dual surfaces.  Every self-dual K\"ahler metric is
weakly selfdual.  These metrics were classified by Bryant in [B].
From [A-C-G-1] it follows that self dual K\"ahler metrics are
orthotoric or of Calabi type and in fact are ambi-K\"ahler. They
are

(a) K\"ahler self-dual metrics of Calabi type over a Riemannian
surface $(\Sigma,g_{\Sigma})$ of constant scalar curvature $k$
where $V(z)=a_1z^4+a_2z^3+kz^2$

(b)  K\"ahler self-dual  metrics of orthotoric type where
$F(x)=lx^3+Ax^2+Bx, G(x)=lx^3+Ax^2+Bx$

(c) complex space forms and a product $\Sigma_c\times\Sigma_{-c}$
of Riemann surfaces of constant scalar curvatures $c$ and $-c$.

\medskip
{\bf Lemma 2.  }  {\it Let    M be a connected QCH K\"ahler
surface which is not Einstein.  Then the following conditions are
equivalent:

(a) The scalar curvature $\tau$ of $\m$ is constant and $
\overline{J}$ is almost K\"ahler

(b) The eigenvalues $\lb,\mu$ of $Ric$ are constant.}

\medskip
{\it Proof.}  (a)$\Rightarrow$(b) Note that
$\rho=\lb\om_1+\mu\om_2$ where $\lb,\mu$ are eigenvalues of $Ric$
and $\om_2=h_J,\om_1=m_J$. Note that $d\om_1+d\om_2=0$ and
$$(\mu-\lb)d\om_1=d\lb\w\om_1+d\mu\w\om_2\tag 2.18$$ Note that $
\overline{J}$ is almost K\"ahler if and only if $d\om_1=0$.  Hence
from (2.7) we get $p_{\De}(\n\lb)=0,p_{\E}(\n\mu)=0$.  Since
$\tau$ is constant we get $\n\lb=-\n\mu$ in an open set
$U=\{x:\lb(x)\ne\mu(x)\}$. Thus $\n\lb=\n\mu=0$ in $U$ and
consequently $U=M$ and $\lb,\mu$ are constant.

$(b)\Rightarrow(a)$ This implication is trivial.$\k$
\medskip
Now we give a classification of locally homogeneous QCH K\"ahler
surfaces.

\medskip

{\bf Proposition 2.} {\it Let $\m$ be a QCH locally homogeneous
manifold. Then the following cases occur:

(a)  $\m$ has constant holomorphic curvature (hence is locally
symmetric and self-dual)

(b) $\m$ is locally a product of two Riemannian surfaces of
constant scalar curvature

(c)  $\m$ is locally isometric to a unique 4-dimensional proper
3-symmetric space.}

{\it Proof.}  If $(M,g)$ is Einstein locally homogeneous
4-manifold then is locally symmetric (see [Jen]). A locally
irreducible locally symmetric K\"ahler surface is self-dual.(see
[D-1]). If $(M,g)$ is not Einstein then using Lemma we see that
$(M,g, \overline{J})$ is an almost K\"ahler manifold satisfying
the Gray condition $G_2$. Hence $||\n \overline{J}||$ is constant
on $M$ and in the case $||\n \overline{J}||\ne0$ it is strictly
almost K\"ahler manifold satisfying $G_2$. Such manifolds are
classified in [A-A-D] and are locally isometric to a proper
3-symmetric space. Note that they are K\"ahler in an opposite
orientation. If $||\n \overline{J}||=0$ then the case (b)
holds.$\k$

{\it Remark.}  A Riemannian 3-symmetric space is a manifold
$(M,g)$ such that for each $x\in M$ there exists an isometry
$\th_x\in Iso(M)$ such that $\th_x^3=Id$ and $x$ is an isolated
fixed point. On a such manifold there is a natural canonical
$g$-ortogonal almost complex structure $\bJ$ such that all $\th_x$
are holomorphic with respect to $\bJ$.  Such structure in
dimension $4$ is almost K\"ahler and satisfies the Gray condition
$G_2$.  The example of 3-symmetric 4-dimensional Riemannian space
with non-itegrable structure $\bJ$ was constructed by O. Kowalski
in [Ko],Th.VI.3. This is the only proper generalized symmetric
space in dimension 4. This example is defined on $\Bbb
R^4=\{x,y,u,v\}$ by the metric
$$\gather g=(-x+\sqrt{x^2+y^2+1})du^2+(x+\sqrt{x^2+y^2+1})dv^2-2ydu\odot
dv\\
+[\frac{(1+y^2)dx^2+(1+x^2)dy^2-2xydx\odot
dy}{1+x^2+y^2}]\endgather$$ It admits a K\"ahler structure $J$ in
an opposite orientation.
\medskip
{\bf Proposition 3.}  {\it Let $\m$ be a QCH K\"ahler surface.  If
$(M,g)$ is conformally Einstein then the almost Hermitian
structure $ \overline{J}$ is Hermitian or $\m$ is self-dual.}
\medskip
{\it Proof.} Let us assume that $(M,g_1)$ is an Einstein manifold
where $g_1=f^2g$. Then $(M,g_1)$ is an Einstein manifold with
degenerate half-Weyl tensor $W^-$. Consequently $W^-=0$ or
$W^-\ne0$ everywhere. In the second case the metric $$
(g_1(W^-,W^-))^{\frac13}g_1$$ is a K\"ahler metric with respect to
$\overline{J}$. Thus $ \overline{J}$ is Hermitian and conformally
K\"ahler.$\k$
\medskip
{\it Remark.}  Every QCH K\"ahler surface is a holomorphically
pseudosymmetric K\"ahler manifold. (see [O],[J-1] ).  In fact from
[J-1] it follows that $R.R=(a+\frac b2)\Pi.R$. Hence in the case
of QCH K\"ahler surfaces we have
$$  R.R=\frac16(\tau-\kappa)\Pi.R\tag 2.19$$
where $\tau$ is the scalar curvature of $\m$ and $\kappa$ is the
conformal scalar curvature of $(M,g, \overline{J})$. Note that
(2.19)  is the obstruction for a K\"ahler surface to have a
negative almost complex $\bJ$ structure satisfying the Gray
condition ($G_2$). In an extremal situation where $(M,g,\bJ)$ is
K\"ahler we have $R.R=0$.

Now we classify QCH K\"ahler surfaces for which $a,b,c$ are all
constant.  Then $\lb,\mu$ are constant and if $(M,g)$ is not
Einstein the almost complex structure $ \overline{J}$ is almost
K\"ahler. Hence $(M,g,\overline{J})$ is a $G_2$ almost K\"ahler
manifold.  Consequently $|\n \overline{\om}|$ is constant and $\m$
is a product of two Riemannian surfaces of constant scalar
curvature or is a proper 3-symmetric space.  If $(M,g)$ is
Einstein then $\kappa=2c$ is constant and
$|W^-|^2=\frac1{24}\kappa^2$ is constant. Thus $\kappa=0$ and $\m$
has constant holomorphic curvature (is a real space form) or by
[D-1] the manifold $(M,g,\overline{J})$ is K\"ahler hence $\m$ is
a product of  two Riemannian surfaces of constant scalar
curvature. Note that for a proper 3-symmetric space we have
$\delta=\frac{\kappa}4$ for the distribution  $\De$ perpendicular
to the K\"ahler nullity of $ \bJ $(see [A-A-D]), thus $
b=2\delta-\frac{\kappa}2=0$ and
$a=\frac16(\tau-\kappa)=-\frac12|\n \overline{\om}|^2$.  Since
$\mu=0$ $c=-\frac32a$ and $\tau=-\kappa$ where $\kappa=\frac32|\n
\overline{\om}|^2$.   Hence
$$R.R=-\frac{\kappa}3\Pi.R\tag2.20$$
where $\kappa=\frac32|\n \overline{\om}|^2$ is constant.
Summarizing we have proved
\medskip
{\bf Proposition 4.}  {\it Let us assume that $\m$ is a QCH
K\"ahler surface with constant $a,b,c$. Then the following cases
occur:

(a)  $\m$ has constant holomorphic curvature (hence is locally
symmetric and self-dual)

(b) $\m$ is locally a product of two Riemannian surfaces of
constant scalar curvature

(c)  $\m$ is locally isometric to a unique 4-dimensional proper
3-symmetric space and $a=-\frac13\kappa,b=0,c=\frac12\kappa$ where
$\kappa=\frac32|\n \overline{\om}|^2$ is constant scalar curvature
of $(M,g,\bJ)$, consequently
$R=-\frac13\kappa\Pi+\frac12\kappa\Psi$.}

\medskip
{\it Remark.}  We consider above the proper 3-symmetric space as a
QCH manifold with respect to the distribution  $\De$ perpendicular
to the K\"ahler nullity of $ \bJ $.  If we consider it as a QCH
manifold with respect to the distribution  $\E=\DE$ then
$R=\frac16\kappa\Pi-\kappa\Phi'+\frac12\kappa\Psi'$ (see Prop.1.).

\bigskip
\centerline{\bf References.}
\par
\medskip
\cite{B}  Bryant R.  {\it Bochner-K\"ahler metrics} J. Amer. Math.
Soc.14 (2001) , 623-715.

\cite{A-C-G-1} V. Apostolov,D.M.J. Calderbank, P. Gauduchon {\it
The geometry of weakly self-dual K\"ahler surfaces}  Compos. Math.
135, 279-322, (2003)
\par
\medskip
\cite{A-C-G-2} V. Apostolov,D.M.J. Calderbank, P. Gauduchon{\it
Ambitoric geometry I: Einstein metrics and extremal ambik\"ahler
structures} arxiv
\par
\medskip
\cite{A-A-D} V. Apostolov, J. Armstrong and T. Draghici {\it Local
ridigity of certain classes Almost K\"ahler 4-manifolds} Ann.
Glob. Anal. and Geom  21; 151-176,(2002)
\par
\medskip
\cite{A-G} V. Apostolov, P. Gauduchon { The Riemannian
Goldberg-Sachs Theorem}  Internat. J. Math. vol.8, No.4,
(1997),421-439
\par
\medskip
\cite{Bes}  A. L. Besse {\it Einstein manifolds}, Ergebnisse,
ser.3, vol. 10, Springer-Verlag, Berlin-Heidelberg-New York, 1987.

\par
\medskip
\cite{D-1} A. Derdzi\'nski, {\it Self-dual K\"ahler manifolds and
Einstein manifolds of dimension four }, Compos. Math.
49,(1983),405-433
 \par
\medskip
\cite{D-2}  A. Derdzi\'nski,  {\it Examples of K\"ahler and
Einstein self-dual metrics on complex plane} Seminar Arthur Besse
1978/79.

\par
\medskip
\cite{G-M-1} G.Ganchev, V. Mihova {\it K\"ahler manifolds of
quasi-constant holomorphic sectional curvatures}, Cent. Eur. J.
Math. 6(1),(2008), 43-75.
\par
\medskip
\cite{G-M-2} G.Ganchev, V. Mihova {\it Warped product K\"ahler
manifolds and Bochner-K\"ahler metrics}, J. Geom. Phys. 58(2008),
803-824.
\par
\medskip
\cite{J-1} W. Jelonek,{ Compact holomorphically pseudosymmetric
K\"ahler manifolds} Coll. Math.117,(2009),No.2,243-249.
\par
\medskip
\cite {J-2} W.Jelonek {\it K\"ahler manifolds with quasi-constant
holomorphic curvature}, Ann. Glob. Anal. and Geom, vol.36, p.
143-159,( 2009)
\par
\medskip
\cite{J-3} W. Jelonek, {\it Holomorphically pseudosymmetric
K\"ahler metrics on $\Bbb{CP}^n$} Coll.
Math.127,(2012),No.1,127-131.
\par
\medskip
\cite{Jen} G.R.Jensen {\it Homogeneous Einstein manifolds of
dimension four} J. Diff. Geom. 3,(1969) 309-349.
\par
\medskip
\cite{Ko} O. Kowalski {\it Generalized symmetric spaces} Lecture
Notes in Math. 805, Springer, New York,1980.
\par
\medskip
\cite{O} Z. Olszak, {\it Bochner flat K\"ahlerian manifolds with a
certain condition on the Ricci tensor} Simon Stevin 63,
(1989),295-303
\par
\medskip
\cite{K-N} S. Kobayashi and K. Nomizu {\it Foundations of
Differential Geometry}, vol.2, Interscience, New York  1963

\par
\medskip

\par
\medskip
Institute of Mathematics

Cracow University of Technology

Warszawska 24

31-155      Krak\'ow,  POLAND.

 E-mail address: wjelon\@pk.edu.pl
\bigskip

\enddocument